\documentclass[12pt]{article}
\usepackage{amsfonts}
\usepackage{amsmath}
\usepackage{a4wide}
\usepackage{amssymb}

\usepackage[latin1]{inputenc}

\newtheorem{lem}{Lemma}[section]
\newtheorem{theo}[lem]{Theorem}
\newtheorem{cor}[lem]{Corollary}
\newtheorem{propo}[lem]{Proposition}
\newtheorem{remar}[lem]{Remark}
\newtheorem{defin}[lem]{Definition}
\newtheorem{conje}[lem]{Conjecture}

\newcommand\proof{\paragraph{Proof}}

\newcommand\F{{\mathbb F}}
\newcommand\Q{{\mathbb Q}}

\newcommand\Z{{\mathbb Z}}

\newcommand\eps\varepsilon

\renewcommand\mod\bmod

\newcommand\gerl{{\mathfrak l}}

\newcommand\gal{{\mathrm {Gal}}}

\title{Note on the diophantine equation $X^t+Y^t=BZ^t${}\footnote{2010 Mathematics Subject Classification:
    Primary 11D41, secondary 11R04,11R18}}

\author{Benjamin Dupuy (Talence)\footnote{Universit\'e Bordeaux
    1, Institut de Math\'ematiques, 351 cours de la Lib\'eration,
    33405 Talence, France}}

\date{}
\begin{document}

\maketitle

\hspace{2cm}

\section{Introduction}\label{sec:1}
Let $t>3$ be a prime number, $B$ be a nonzero rational integer. Consider the equation 
\begin{equation}
\label{eq:1}
X^t+Y^t=BZ^t
\end{equation}
where $X,$ $Y,$ $Z$ are coprime non zero rational
integers.
\begin{defin}
Let $t>3$ be a prime number. We say that $t$ is a good prime number if
and only if 

\medskip

\begin{itemize}
\item its index irregularity $\iota(t)$ is equal to zero 


or


\item $t\nmid{h_{t}^+}$ and none of the
Bernoulli numbers $B_{2nt}$, $n=1,\ldots{,\frac{t-3}{2}}$ is divisible
by $t^3$.
\end{itemize}
\end{defin}

\medskip

For a prime number $t$ with $t<12.10^6$, it has been recently proved that none of the
Bernoulli numbers $B_{2nt}$, $n=1,\ldots{,\frac{t-3}{2}}$ is divisible
by $t^3$ (see \cite{Buhl}). Furthermore, $h_{t}^+$ is prime to $t$ for
$t<7.10^6$.  In particular, every prime number $t<7.10^6$ is a
good prime number in the previous meaning. 

\medskip

As usual, we denote by $\phi$ the Euler's function. \textbf{For the following,
we fix $t>3$ a good prime number}, and a rational integer $B$ prime to
$t$, such that for every prime number $l$
dividing $B$, we have $-1\mod t\in{<l\mod t>}$ the subgroup of
$\F_{t}^{\times}$ generated by $l \mod t$. For example, it is the case
if for every prime number $l$ dividing $B$, $l\mod t$ is not a square.

 In this paper, using
a descent method on the number of prime ideals, we prove the
following theorem
\begin{theo}
\label{th1} The equation (\ref{eq:1}) has no solution in pairwise relatively prime non
zero integers $X,$ $Y,$ $Z$ with $t|Z$. 
\end{theo}

\medskip

In particular, using a recent result of Bennett et al, we deduce the
\medskip
\begin{cor} 
\label{corth1} Suppose that $B^{t-1}\neq 2^{t-1}\mod t^2$ and $B$ has a divisor
  $r$ such that $r^{t-1}\neq 1\mod t^2$. Then the equation
  (\ref{eq:1}) has no solution in pairwise relatively prime non zero integers $X,$ $Y,$ $Z$. 
\end{cor}

\section{Proof of the theorem}

First, we suppose that $\iota(t)=0$. Let us prove the following lemma.

\begin{lem} \label{des}
Let $\zeta$ be a primitive $t$-th root of unity and $\lambda=(1-\zeta)(1-\overline{\zeta})$. Suppose there exist algebraic integers $x,y,z$ in the ring $\Z[\zeta+\overline{\zeta}]$, an integer $m\geq{t}$, and a unit $\eta$ in $\Z[\zeta+\overline{\zeta}]$ such that $x$, $y$, $z$ and $\lambda$ are pairwise coprime and verify
\begin{equation}
\label{eqKummer} x^t+y^t=\eta\lambda^{m}Bz^t.
\end{equation}
Then $z$
is not a unit of $\Z[\zeta+\overline{\zeta}]$. Moreover, there
exist algebraic integers $x',y',z'$ in $\Z[\zeta+\overline{\zeta}]$, an integer $m'\geq{t}$, and a unit $\eta'$ in $\Z[\zeta+\overline{\zeta}]$ such that $x'$, $y'$, $z'$, $\lambda$ and $\eta'$ verify the same properties. The algebraic number $z'$ divides $z$ in $\Z[\zeta]$. The number of prime ideals of $\Z[\zeta]$ counted with multiplicity and dividing $z'$ is strictly less than that of $z$.
\end{lem}

\proof
The equation (\ref{eqKummer}) becomes
\begin{align*}
\left(x+y\right)\prod_{a=1}^{t-1}{\left(x+\zeta^a
y\right)}=\eta\lambda^{m}Bz^t.
\end{align*}
By hypothesis, for every prime number $l$ dividing $B$, we have $-1\mod
t\in{<l\mod t>}$. In particular $B$ is prime to
$\frac{x^t+y^t}{x+y}$. In fact, suppose there exist $\gerl$ a prime
factor of $B$ in $\Z[\zeta]$ such that $\gerl
|\frac{x^t+y^t}{x+y}$. Then there exist $a\in{\{1,\ldots{,t-1}\}}$,
such that $\gerl |x+\zeta^a y$. Let $l$ be the rational prime number
under $\gerl$. Since $-1\mod t$ is an element of the subgroup of
$\F_{t}^{\times}$ generated by $l\mod t$, we deduce that the decomposition group of
$\gerl$ contains the complex conjugation
$j\in{\gal{\left(\Q(\zeta)/\Q\right)}}$ that is $\gerl^{j}=\gerl$. Particularly, $\gerl
|x+\zeta^a y$ implies that $\gerl|x+\zeta^{-a}y$ since $x$, $y$
are real. So $\gerl|(\zeta^a -\zeta^{-a})y$. Since $\gerl$ is a prime
ideal, we deduce that $\gerl|y$ or $\gerl|\zeta^a -\zeta^{-a}$. But since $x$ and $y$ are
coprime, $y$ is prime to $\gerl$. Since $(B,p)=1$ and $\zeta^a-\zeta^{-a}$ is a
generator of the only prime ideal of $\Z[\zeta]$ above $p$, we can not
have $\gerl|\zeta^a -\zeta^{-a}$: we get a contradiction. So $B$ and
$\frac{x^t+y^t}{x+y}$ are coprime as claimed. In fact, we have proved
the following result: $B$ is prime to every factor of the form
$\frac{a^t+b^t}{a+b}$ where $a$ and $b$ are coprime elements of $\Z[\zeta+\overline{\zeta}]$.

Then $B|x+y$ in $\Z[\zeta]$. Therefore we get
\begin{align*}
\frac{x+y}{B}\prod_{a=1}^{t-1}{\left(x+\zeta^a
y\right)}=\eta\lambda^{m}z^t.
\end{align*}
Following the same method\footnote{Recall that $t\nmid{h_{t}^+}$ since
$\iota(t)=0$} as in section $9.1$ of \cite{wa}, one can show that there exist real units
$\eta_{0},\eta_{1},\ldots{,\eta_{t-1}}\in{\Z[\zeta+\overline{\zeta}]^{\times}}$
and algebraic integers $\rho_{0}\in{\Z[\zeta+\overline{\zeta}]}$,
$\rho_{1},\ldots{,\rho_{p-1}}\in{\Z[\zeta]}$ such that
\begin{equation}
\label{relkum}
x+y=\eta_{0}B\lambda^{m-\frac{t-1}{2}}\rho_{0}^t,\quad
\frac{x+\zeta^a y}{1-\zeta^a}=\eta_{a}\rho_{a}^t,\quad
a=1,\ldots{,t-1.}
\end{equation}
Let us show that $z$ is not a unit. As $\rho_{1}$ divides $z$ in $\Z[\zeta]$, it is thus enough to show that $\rho_{1}$ is not one. Put $\alpha= \frac{x+\zeta y}{1-\zeta}$. One has
\begin{align*}
\alpha=-y+\frac{x+y}{1-\zeta}\equiv -y\mod(1-\zeta)^2.
\end{align*}
So $\frac{\overline{\alpha}}{\alpha}\equiv 1\mod
(1-\zeta)^2$. Suppose that $\rho_ {1}$ is a unit. Then, the quotient
$\frac{\overline{\rho_{1}}^t}{\rho_{1}^t}$ is a unit of modulus
$1$ of the ring $\Z[\zeta]$, thus a root of the unity of this ring by Kronecker theorem. However, the only roots of the unity of $ \Z [\zeta]$ are the $2t$-th roots of the unity (see \cite{wa}). As the unit $\eta_{1}$ is real, thus there exists an integer $l$ and $\epsilon= \pm 1$ such as
$\frac{\overline{\eta_{1}}\cdot\overline{\rho_{1}}^t}{\eta_{1}\cdot\rho_{1}^t}=\frac{\overline{\rho_{1}}^t}{\rho_{1}^t}=\epsilon\zeta^l$.
Therefore, we have
\begin{align*}
\frac{\overline{\alpha}}{\alpha}=\epsilon\zeta^l.
\end{align*}
As $\frac{\overline{\alpha}}{\alpha}\equiv 1\mod (1-\zeta)^2$,
we get $\epsilon\zeta^l\equiv 1\mod (1-\zeta)^2$, so
$\epsilon\zeta^l=1$, i.e. $\frac{\overline{\alpha}}{\alpha}=1$. So
\begin{align*}
\frac{x+\zeta
y}{1-\zeta}=\frac{x+\overline{\zeta}y}{1-\overline{\zeta}},
\end{align*}
because $x$ and $y$ are real numbers. From this equation, we deduce that \begin{align*}
\frac{x+\zeta y}{1-\zeta}=\frac{\zeta x+y}{\zeta-1},\quad i.e.\quad
(x+y)(\zeta+1)=0.
\end{align*}
We get a contradiction. So the algebraic integer $\rho_{1}$ (and then $z$) is not a unit. This completes the proof of first part of the lemma.

Let us prove the existence of $x'$, $y'$, $z'$, $\eta',$ and $m'$. It
is just an adaptation of the computations done in paragraph $9.1$ of
chapter $9$ of \cite{wa} for the second case of the Fermat
equation. We give here the main ideas. Let
$a\in{\left\{1,\ldots{,p-1}\right\}}$ be a fixed integer. We take
$\lambda_{a}=(1-\zeta^a)(1-\zeta^{-a})$. By \eqref{relkum}, there
exist a real unit $\eta_{a}$ and $\rho_{a}\in{\Z[\zeta]}$ such that
\begin{align*}
\frac{x+\zeta^a y}{1-\zeta^a}=\eta_{a}\rho_{a}^t,
\end{align*}
and taking the conjugates (we know that $x,y\in{\mathbb{R}}$), we have
\begin{align*}
\frac{x+\zeta^{-a} y}{1-\zeta^{-a}}=\eta_{a}\overline{\rho_{a}}^t.
\end{align*}
Thus
$$
x+\zeta^a y=(1-\zeta^a)\eta_{a}\rho_{a}^t,\;\;
x+\zeta^{-a} y=(1-\zeta^{-a})\eta_{a}\overline{\rho_{a}}^t.
$$
Multiplying the previous equalities, we obtain
\begin{equation}
\label{eqkum2}
x^2+y^2+\left(\zeta^a+\zeta^{-a}\right)xy=\lambda_{a}\eta_{a}^2\left(\rho_{a}\overline{\rho_{a}}\right)^t.
\end{equation}
Taking the square of
$x+y=\eta_{0}B\lambda^{m-\frac{t-1}{2}}\rho_{0}^t$ gives
\begin{equation}
\label{eqkum3}
x^2+y^2+2xy=\eta_{0}^2B^2\lambda^{2m-t+1}\rho_{0}^{2t}.
\end{equation}
The difference between equations \eqref{eqkum3}, \eqref{eqkum2} and then the division by $\lambda_{a}$ give
\begin{equation}
\label{eqkum4}
-xy=\eta_{a}^2\left(\rho_{a}\overline{\rho_{a}}\right)^t-\eta_{0}^2B^2\lambda^{2m-t+1}\rho_{0}^{2t}\lambda_{a}^{-1}.
\end{equation}
As $t>3$, there exists an integer $b\in{\{1,\ldots{,t-1}\}}$ such that $b\neq \pm a\mod t$. For this integer $b$, we get
\begin{equation}
\label{eqkum5}
-xy=\eta_{b}^2\left(\rho_{b}\overline{\rho_{b}}\right)^t-\eta_{0}^2B^2\lambda^{2m-t+1}\rho_{0}^{2t}\lambda_{b}^{-1}.
\end{equation}
The difference of between equations \eqref{eqkum4} and \eqref{eqkum5} gives after simplifying
\begin{align*}
\eta_{a}^2\left(\rho_{a}\overline{\rho_{a}}\right)^t-\eta_{b}^2\left(\rho_{b}\overline{\rho_{b}}\right)^t=\eta_{0}^2B^2\lambda^{2m-t+1}\rho_{0}^{2t}\left(\lambda_{a}^{-1}-\lambda_{b}^{-1}\right).
\end{align*}
But as $b\neq\pm a\mod t$, we have
$\lambda_{a}^{-1}-\lambda_{b}^{-1}=\frac{\left(\zeta^{-b}-\zeta^{-a}\right)\left(\zeta^{a+b}-1\right)}{\lambda_{a}\lambda_{b}}=\frac{\delta'}{\lambda}$,
where $\delta'$ is a unit. We know that $\lambda_{a}$, $\lambda_{b}$ and
$\lambda$ are real numbers, then the unit $\delta'$ is a real unit. So there exists a real unit $\eta'=\frac{\delta'\cdot\eta_{0}^2}{\eta_{b}^2}$ such that
\begin{equation}
\label{eqkum6}
\left(\frac{\eta_{a}}{\eta_{b}}\right)^2\left(\rho_{a}\overline{\rho_{a}}\right)^t
+\left(-\rho_{b}\overline{\rho_{b}}\right)^t=\eta'B^2\lambda^{2m-t}\left(\rho_{0}^2\right)^t.
\end{equation}
The condition $\iota(t)=0$ implies that $\frac{\eta_{a}}{\eta_{b}}$
is a $t$-th power in $\Z[\zeta+\overline{\zeta}]$. Thus there exists $\xi\in{\Z[\zeta+\overline{\zeta}]}$ such that
$\frac{\eta_{a}}{\eta_{b}}=\xi^t$. In fact, we know that
\begin{align*}
\eta_{a}\rho_{a}^t=\frac{x+\zeta^a y}{1-\zeta^a},\quad
x+y=\eta_{0}B\lambda^{m-\frac{t-1}{2}}\rho_{0}^t\equiv 0\mod
(1-\zeta)^{2m-t+1}.
\end{align*}
Then
\begin{align*}
\eta_{a}\rho_{a}^t=-y+\frac{x+y}{1-\zeta^a}\equiv-y\mod(1-\zeta)^{2m-t}
\equiv-y\mod t.
\end{align*}
Also $\eta_{b}\rho_{b}^t\equiv-y\mod t$, where
$\frac{\eta_{a}}{\eta_{b}}\equiv\left(\frac{\rho_{b}}{\rho_{a}}\right)^t\mod
t$. But Lemma $1.8$ in \cite{wa} shows that
$\left(\frac{\rho_{b}}{\rho_{a}}\right)^t$ is congruent to an integer $l\in{\Z}$ modulo $t$, therefore the existence of an integer $l$
such that
\begin{align*}
\frac{\eta_{a}}{\eta_{b}}\equiv l\mod t,\quad l\in{\Z}.
\end{align*}
By Theorem $5.36$ of \cite{wa}, the unit $\frac{\eta_{a}}{\eta_{b}}$ is a $t$-th power in $\Z[\zeta]$ so we have the existence of $\xi_{1}\in{\Z[\zeta]}$ such that $\frac{\eta_{a}}{\eta_{b}}=\xi_{1}^t$. As the unit $\frac{\eta_{a}}{\eta_{b}}$ is real, one has
\begin{align*}
\xi_{1}^t=\overline{\xi_{1}}^t.
\end{align*}
Therefore, there exists an integer $g$ such that
$\overline{\xi_{1}}=\zeta^g\xi_{1}$. Taking
$\xi=\zeta^{gh}\xi_{1}$ where $h$ is the inverse of $2\mod t$, we have
\begin{align*}
\overline{\xi}=\xi,\quad \xi^t=\xi_{1}^t=\frac{\eta_{a}}{\eta_{b}},
\end{align*}
i.e. $\frac{\eta_{a}}{\eta_{b}}=\xi^t$, where
$\xi\in{\Z[\zeta+\overline{\zeta}]}$. We put
\begin{align*}
x'=\xi^2 \rho_{a}\overline{\rho_{a}},\quad
y'=-\rho_{b}\overline{\rho_{b}},\quad z'=\rho_{0}^2,\quad m'=2m-t.
\end{align*}
One can verify that
\begin{align*}
x'^t+y'^t=\eta' B^2\lambda^{m'}z'^t.
\end{align*}
Obviously, $B^2$ is prime to $t$ and  for all prime $l$
dividing $B^2$, we have $-1\mod t\in{<l\mod t>}$ the subgroup of
$\F_{t}^{\times}$ generated by $l \mod t$.  Moreover, one have already seen that the algebraic integer $\rho_{1}$ is not a unit in
$\Z[\zeta]$. As $\rho_{0}\rho_{1}$ divides $z$ in $\Z[\zeta]$,
the number of prime ideals counted with multiplicity and dividing
$z'$ in $\Z[\zeta]$ is then strictly less than that of $z$ and $m'=2m-t\geq{2t-t=t}$. This completes the proof of the lemma.$\square$

Now let $(X,Y,Z)$ be a solution of (\ref{eq:1}) in pairwise relatively prime non
zero integers with $t|Z$. Let $Z=t^{v}Z_{1}$ with $t\nmid{Z_{1}}$. Equation (\ref{eq:1}) becomes
\begin{align*}
X^t+Y^t=Bt^{tv}Z_{1}^t.
\end{align*}
Let $\zeta$ be a primitive $t$-th root of unity and $\lambda=(1-\zeta)(1-\overline{\zeta})$. The previous equation
becomes
\begin{align*}
X^t+Y^t=B\frac{t^{tv}}{\lambda^{tv\frac{t-1}{2}}}\lambda^{tv\frac{t-1}{2}}Z_{1}^t.
\end{align*}
The quotient $\eta=\frac{t^{tv}}{\lambda^{tv\frac{t-1}{2}}}$ is a real
unit in the ring $\Z[\zeta+\overline{\zeta}]$. Take
$m=tv\frac{t-1}{2}\geq{t}$. We have just proved that there exist
$\eta\in{\Z[\zeta+\overline{\zeta}]^{\times}}$ and an integer $m\geq{t}$ such that
\begin{equation}
\label{eqkumdep} X^t+Y^t=\eta B\lambda^{m}Z_{1}^t,
\end{equation}
where $X$, $Y$, $\lambda$ and $Z_{1}$ are coprime.

We can apply the lemma (\ref{des}) to equation \eqref{eqkumdep}. By induction, one can prove
the existence of the sequence of algebraic $Z_{i}$ such that
$Z_{i+1}|Z_{i}$ in $\Z[\zeta]$ and the number of prime factors in
$\Z[\zeta]$ is strictly decreasing. So there is $n$ such that $Z_{n}$
is a  unit. But Lemma \ref{des} indicates that  each of the $Z_{i}$ is
not a unit. We get a contradiction. The theorem is proved in the case
$\iota(t)=0$. 

In the other case, $(t,h_{t}^+)=1$ and none of the
Bernoulli numbers $B_{2nt}$, $n=1,\ldots{,\frac{t-3}{2}}$ is divisible
by $t^3$. In particular, with the notation of the proof of the lemma,
it exists $\xi\in{\Z[\zeta+\overline{\zeta}]}$ such that
$\frac{\eta_{a}}{\eta_{b}}=\xi^t$ (see \cite{wa}, page $174-176$). So the results of the previous lemma
are valid in the second case. We conclude as before. The theorem is proved.

\section{Proof of the corollary}

Let $X,$ $Y,$ $Z$ be a solution in pairwise relatively prime non zero
integers of the equation (\ref{eq:1}). By the theorem, the integer $Z$
is prime to $t$. Furthermore, $B\phi(B)$ is coprime to $t$, $B^{t-1}\neq
2^{t-1}\mod t^2$ and $B$ has a divisor $r$ such that $r^{t-1}\neq 1\mod
t^2$. So by the theorem $4.1$ of \cite{Be}, the equation (\ref{eq:1})
has no solution for such $t$ and $B$.

\paragraph{Acknowledgments} I thank Professor Yuri Bilu and Florian
Luca for very helpful suggestions.

\end{document}